\documentclass[12pt,a4paper]{article}
\usepackage{amsmath, amstext, amsxtra, amsthm, amscd, amsgen, amsbsy, amsopn,amsfonts, latexsym, amssymb, amscd, array}

\newtheorem{theorem}{Theorem}[section]
\newtheorem{proposition}[theorem]{Proposition}
\newtheorem{lemma}[theorem]{Lemma}
\newtheorem{rrremark}[theorem]{Remark}
\newtheorem{finalremark}[theorem]{Final remark}

\newtheorem{corollary}[theorem]{Corollary}

\newtheorem{observation}[theorem]{Observation}

\begin{document}

\def\Z{{\mathbb Z}}
\def\Zp#1{{\mathbb Z}/#1{\mathbb Z}}
\def\R{{\mathbb R}}
\def\N{{\mathbb N}}
\def\H{{\mathbb H}}
\def\Q{{\mathbb Q}}
\def\C{{\mathbb C}}

\def\proof#1{\bigskip
\noindent {\bf Proof:} #1

\bigskip
\noindent}

\def\proofend{\hbox to 1em{\hss}\hfill $\blacksquare $}

\title{On the topology of scalar-flat manifolds}
\author{Anand Dessai}
\date{}
\maketitle
\noindent
\abstract{\noindent Let $M$ be a simply-connected closed manifold of
dimension $\geq 5$ which does not admit a metric with positive scalar
curvature. We give necessary conditions for $M$ to admit a scalar-flat
metric. These conditions involve the first Pontrjagin class and the
cohomology ring of $M$. As a consequence any simply-connected
scalar-flat manifold of dimension $\geq 5$ with vanishing first
Pontrjagin class admits a metric with positive scalar curvature. We 
also describe some relations between scalar-flat metrics, almost 
complex structures and the free loop space. 

\bigskip
\begin{center}1. {\it Introduction}\end{center}
\setcounter{section}{1}
In this paper we give restrictions on the possible topological type of
connected closed scalar-flat manifolds. A Riemannian manifold is
called scalar-flat if its scalar curvature vanishes identically. One
way to obtain such manifolds (in dimension $\geq 3 $) is to start with
a metric with positive scalar curvature and decrease the scalar
curvature by changing the metric. By the results of Kazdan-Warner (cf.
\cite{KaWa75}) one can deform the metric globally by a conformal
change to obtain a scalar-flat metric. Recently Lohkamp (cf.
\cite{Lo97}) has shown that scalar decreasing deformations can even be
carried out locally to achieve changes of the scalar curvature
arbitrarily close to a prescribed decrease.

However, not every scalar-flat metric arises from a metric with 
positive scalar curvature. The positive energy theorem and the 
Lichnerowicz formula for $Spin$-manifolds imply that it is in general 
impossible to increase scalar curvature locally or globally. In the 
following we shall call a scalar-flat manifold {\it strongly 
scalar-flat} if it does not admit a metric with positive scalar 
curvature. By a result of Bourguignon such manifolds are Ricci-flat. 
Otherwise one could deform the metric using the Ricci tensor to obtain 
a new metric which would be pointwise conformal to one with positive 
scalar curvature (cf. \cite{KaWa75}). 

Although strongly scalar-flat manifolds seem to be quite special very
little is known about their to\-po\-logical type. In contrast the
surgery theorem of Gromov-Lawson and Schoen-Yau (cf. \cite{GrLa80},
\cite{ScYa79}) and Stolz' solution of the Gromow-Lawson conjecture
(cf. \cite{St90}) give a complete topological classification for
simply-connected manifolds of dimension $\geq 5$ which admit a metric
with positive scalar curvature.

In Theorem \ref{main} we give partial information on the topological
type of strongly scalar-flat manifolds. Our result shows that many
manifolds are not strongly scalar-flat. In particular, it implies the
following
\begin{theorem}\label{firstpont} Let $M$ be a simply-connected scalar-flat manifold of dimension
$\geq 5$. If the first real Pontrjagin class of $M$ vanishes then $M$
is not strongly scalar-flat, i.e. $M$ admits a metric with positive
scalar curvature.
\end{theorem}

\noindent
As a consequence of Stolz' theorem (cf. \cite{St90}) a
simply-connected strongly scalar-flat manifold of dimension $\geq 5$
admits a non-vanishing harmonic spinor. Futaki proved in \cite{Fu93}
that this implies that the manifold decomposes as a Riemannian product
of irreducible strongly scalar-flat manifolds which are Ricci-flat
K\"ahler or have holonomy $Spin(7)$ (cf. also \cite{Hi74}). Theorem
\ref{main} follows from properties of these factors. Futaki used this
decomposition to give an upper bound for the absolute value of the
$\hat A$-genus of strongly scalar-flat manifolds.

The paper is structured in the following way. In Section 2 we state 
our main result on the topology of strongly scalar-flat manifolds (see 
Theorem \ref{main}) and give a partial answer to the question whether 
such manifolds are almost complex (see Proposition \ref{su}). In 
Section 3 we prove these statements and observe that irreducible 
Ricci-flat manifolds with vanishing first Pontrjagin class have 
generic holonomy. In Section 4 we rephrase some of the results in 
terms of the free loop space of the manifold and point out more recent 
results related to this paper.

\bigskip
\begin{center}2. {\it Topology of strongly scalar-flat manifolds}\end{center}
\setcounter{section}{2}
\setcounter{theorem}{0}
In this section we state some necessary conditions for a manifold to
be strongly scalar-flat. We also give a partial answer to the question
whether strongly scalar-flat manifolds are almost complex. In the
following all manifolds shall be smooth, connected and closed.
\begin{theorem}\label{main} Let $M$ be a simply-connected manifold of
dimension $\geq 5$. If $M$ is strongly scalar-flat then $M$ carries a
closed $2$-form $\omega $ and a closed $4$-form $\Omega $ such that
$$\langle e^{\lbrack \omega \rbrack }\cdot e^{\lbrack \Omega \rbrack }\cdot p_1(M),\mu
_M\rangle \neq 0.$$
\end{theorem}

\bigskip
\noindent
Here $\lbrack \omega \rbrack $ and $\lbrack \Omega \rbrack $ denote
the induced real cohomology classes, $\mu _M$ is the fundamental cycle
of $M$ and $\langle \quad ,\quad
\rangle $ denotes the
Kronecker pairing. Theorem \ref{firstpont} follows directly
from Theorem \ref{main}.

We deduce Theorem \ref{main} in the next section from properties of
Ricci-flat K\"ahler manifolds and manifolds with holonomy $Spin(7)$.
First compact examples of $Spin(7)$-holonomy manifolds were
constructed by Joyce (cf. \cite{JoII96}). Very little is known about
their topological type. In all known examples the signature is
divisible by $8$ (cf. \cite{JoII96}, \cite{Jo99}). In particular their
Euler characteristic is even. If this were true for all
$Spin(7)$-holonomy manifolds then any strongly scalar-flat manifold
would be almost complex. A more precise statement is given in the
following
\begin{proposition}\label{su} Let $M$ be a simply-connected strongly scalar-flat manifold of
dimension $\geq 5$. Assume for some scalar-flat metric on $M$ that any
$Spin(7)$-holonomy factor of $M$ has even Euler characteristic. Then
$M$ admits an almost complex structure with vanishing first Chern
class. If $\dim M=4k$ then $sign (M)\equiv (-1)^k\cdot \chi (M) \bmod
4$.
\end{proposition}

\bigskip
\noindent
We close this section with a few remarks.
\begin{rrremark}\label{remark}
\begin{enumerate}
\item The method used to prove Theorem \ref{main} also shows that the odd
 Stiefel-Whitney classes of a simply-connected strong\-ly scalar-flat manifold vanish.
\item In contrast to the case of positive scalar curvature strongly scalar-flat structures are not preserved under surgery. 
In fact, performing a surgery (of any codimension) turns a strongly scalar-flat manifold
with infinite fundamental group into a manifold for which the scalar
curvature of any metric is negative somewhere. The situation for
finite fundamental groups is similar.
\end{enumerate}
\end{rrremark}

\bigskip
\begin{center}3. {\it Holonomy of strongly scalar-flat manifolds}\end{center}
\setcounter{section}{3}
\setcounter{theorem}{0}
In this section we prove Theorem \ref{main}, Proposition \ref{su} and
the following observation concerning the question whether there exist
irreducible Ricci-flat manifolds with generic holonomy (cf. 
\cite{Be87}, p. 19). 
\begin{observation}\label{holonomy} A simply-connected irreducible Ricci-flat manifold with vanishing real first
Pontrjagin class has generic holonomy. 
\end{observation}

\bigskip
\noindent
 Theorem
\ref{main} relies on the following result of Futaki which he deduces
from the surgery theorem (cf. \cite{GrLa80}) and Stolz' solution of
the Gromov-Lawson conjecture (cf. \cite{St90}) using the de Rham
splitting theorem and Berger's classification of holonomy groups (cf.
\cite{Be87}, Chapter 10).
\begin{theorem}\label{futaki}{\rm \bf (\cite{Fu93}, Th. 1)} Let $M$ be a simply-connected strongly scalar-flat manifold of
dimension $\geq 5$. Then $M$ is a $Spin$-manifold and splits as a
Cartesian product of irreducible manifolds, where each factor admits a 
Ricci-flat K\"ahler metric or a Riemannian metric with holonomy 
$Spin(7)$.\proofend 
\end{theorem}

\bigskip
\noindent
In the next proposition we collect those topological properties of
Ricci-flat K\"ahler manifolds and $Spin(7)$-holonomy manifolds which
are used in the proof of Theorem \ref{main}. Note that a holonomy
reduction leads to a corresponding topological reduction of the frame
bundle. We always fix the orientation of the manifold induced by this
reduction.
\begin{proposition}\label{p1nonzero} Let $(M,g)$ be a simply-connected irreducible Riemannian
manifold of dimension $n>0$. Let $H$ be the holonomy group of
$(M,g)$.\begin{enumerate}
\item If $H\subset SU(n/2)$ (i.e. $(M,g)$ is Ricci-flat K\"ahler) and $\omega $ denotes the K\"ahler form
then
$$\langle e^{\lbrack \omega \rbrack }\cdot p_1(M),\mu
_M\rangle <0\quad \text{ and }\quad \langle e^{\lbrack \omega \rbrack },\mu
_M\rangle >0.$$
\item If $H=Spin(7)$ (implying $n=8$) then $M$ admits a closed $4$-form $\Omega $ such that
$$\langle e^{\lbrack \Omega \rbrack }\cdot
p_1(M),\mu _M\rangle < 0\quad \text{ and }\quad \langle e^{\lbrack
\Omega \rbrack },\mu
_M\rangle >0.$$
\end{enumerate}
\end{proposition}

\noindent
{\bf Proof:} Let $H\subset SU(m)$ where $m=n/2$. We note that $\langle
e^{\lbrack \omega \rbrack },\mu
_M\rangle >0$ since the $m$-fold wedge of the K\"ahler form $\omega $ is
 (up to a positive constant) the volume form. Also $c_1(M)=0$ implies 
 that $p_1(M)=-2\cdot c_2(M)$.
Thus $\langle e^{\lbrack \omega \rbrack } \cdot p_1(M),\mu
_M\rangle $ is equal to $\langle \lbrack \omega \rbrack ^{m-2}\cdot
c_2(M),\mu _M\rangle $ times a negative constant. We now use the
following formula of Apte (cf. \cite{Be87}, formula (2.80a))
$$\langle \lbrack \omega \rbrack ^{m-2}\cdot
c_2(M),\mu _M\rangle =C\cdot \int _M\vert R\vert ^2\; \mu _g.$$ Here 
$C$ is a positive constant and $R$ denotes the curvature tensor viewed 
as a symmetric endomorphism on $2$-forms. Note that $(M,g)$ as a 
simply-connected manifold of positive dimension cannot be flat. Hence 
$$\langle e^{\lbrack \omega \rbrack } \cdot p_1(M),\mu
_M\rangle <0.$$
Next let $\dim M=8$ and $H=Spin(7)\hookrightarrow SO(8)$. We adapt the
argument for $G_2$-manifolds given in \cite{JoI96}, p. 333. First we
recall some standard facts about $Spin(7)$-holonomy manifolds (cf.
\cite{Br87}, \cite{JoII96}). Associated to the holonomy reduction is a
parallel closed self-dual $4$-form $\tilde \Omega $ (for an explicit
local description see for example \cite{JoII96}, p. 510). Its
stabilizer at each point is isomorphic to $Spin(7)$. Up to a positive
constant $\tilde \Omega \wedge \tilde \Omega $ is the volume form of
$M$.

The space of $k$-forms $\Gamma (\Lambda ^k(M))$ splits orthogonally
into components para\-metrized by the irreducible representations of
the $Spin(7)$-module $\Lambda ^k(\R ^8)$. In particular $\Gamma
(\Lambda ^2(M))$ decomposes as $\Gamma (\Lambda ^2_{7})\oplus \Gamma
(\Lambda ^2_{21})$, where $\Gamma (\Lambda ^2_{7})$ corresponds to the
$Spin(7)$-representation $\Lambda ^1(\R ^7)$ and $\Gamma (\Lambda
^2_{21})$ corresponds to the adjoint representation of $Spin(7)$.
Local computations (cf. \cite{Br87}, p. 546) show that $\int _M\tilde
\Omega \wedge
\xi
\wedge \xi =-\int _M\vert \xi \vert ^2\; \mu _g$ for any $\xi \in \Gamma (\Lambda ^2_{21})$. So for $\Omega :=-\tilde \Omega $ we get
$$\int _M \Omega \wedge \Omega >0\quad \text{ and }\quad \int _M\Omega \wedge
\xi
\wedge \xi =\int _M\vert \xi \vert ^2\; \mu _g.$$
By the Ambrose-Singer theorem (cf. \cite{Be87}, p. 291) the components
$R_{ij}$ of the curvature tensor $R$ are in $\Gamma (\Lambda
^2_{21})$. Since the first Pontrjagin form $p_1(M,g)$ is given by a
negative multiple of the trace $tr(R\wedge R)=\sum R_{ij}\wedge
R_{ij}$ we conclude
$$\int _M \Omega \wedge p_1(M,g) =C\cdot \int _M\vert R\vert ^2
\; \mu _g,$$
where $C<0$ is a constant. Again the integral is non-zero since
$(M,g)$ is not flat. Hence, $\langle e^{\lbrack \Omega \rbrack }\cdot p_1(M),\mu _M \rangle =\int _M \Omega
\wedge p_1(M,g)< 0$.\proofend

\bigskip
\noindent
{\bf Proof of Theorem \ref{main}:} Let $M=K_1\times \ldots \times
K_r\times J_1\times \ldots \times J_s$ be the splitting of the strongly scalar-flat manifold $M$ into irreducible factors, where the holonomy
of $K_i$ is contained in the special unitary group (i.e. $K_i$ is
Ricci-flat K\"ahler) and the holonomy of $J_j$ is $Spin(7)$ (see
Theorem \ref{futaki}). On each factor we choose the orientation
induced by the holonomy reduction. This gives an orientation of $M$.
Let $\omega _i$ denote the K\"ahler form of $K_i$ and let $\Omega _j$
denote the $4$-form given in Proposition \ref{p1nonzero}. We put $\omega :=\sum _i\omega _i$, $\Omega :=\sum _j\Omega _j$, $A_i:=\langle e^{\lbrack \omega _i\rbrack },\mu _{K_i}\rangle $ ,
$B_j:=\langle e^{\lbrack \Omega _j\rbrack },\mu _{J_j}\rangle $,
$A:=\prod _{i=1}^rA_i$ and $B:=\prod _{j=1}^sB_j$. Next we note that
$\langle e^{\lbrack \omega \rbrack }\cdot e^{\lbrack \Omega \rbrack
}\cdot p_1(M),\mu _M\rangle $ may be written as a sum, where each
summand has one of the following forms
$$\langle e^{\lbrack \omega _i\rbrack }\cdot
p_1(K_i),\mu _{K_i}\rangle \cdot (A\cdot B)/A_i\quad \text{ or }\quad
\langle e^{\lbrack \Omega _j\rbrack }\cdot p_1(J_j),\mu _{J_j}\rangle
\cdot (A\cdot B)/B_j.$$ By Proposition \ref{p1nonzero} each summand is
negative.\proofend

\bigskip
\noindent
{\bf Proof of Observation \ref{holonomy}:} Assume $M$ is an
irreducible Ricci-flat manifold with reduced holonomy. We apply
Berger's classification (cf. \cite{Be87}, Chapter 10) and conclude
that the holonomy of $M$ is $G_2$, $Spin(7)$ or contained in the
special unitary group. In all these cases the first real Pontrjagin
class of $M$ is non-zero giving the desired contradiction (see
\cite{JoI96} and Proposition \ref{p1nonzero}).\proofend
 
\bigskip
\noindent
For the proof of Proposition \ref{su} we use the following
\begin{lemma}\label{acs} Let $N$ be an $8$-dimensional $Spin$-manifold which admits a
topological $Spin(7)$-reduction. Then $N$ admits an almost complex
structure with vanishing first Chern class if and only if the Euler
characteristic $\chi (N)$ is even.
\end{lemma}

\noindent {\bf Proof:} In \cite{He70} Heaps showed, using results of Massey, that an oriented $8$-dimensional manifold $X$ admits a stable almost complex
structure (s.a.c.s.) if and only if $w_2(X)$ is integral and
$w_8(X)\in Sq^2(H^6(X;\Z ))$. Note that the Steenrod operation
$Sq^2:H^6(N ;\Z /2\Z )\to H^8(N ;\Z /2\Z )$ vanishes identically since
it is given by multiplication with $w_2(N)=0$. Since $\chi (N)$
reduces to $w_8(N)$ modulo $2$ we conclude that $N$ admits a s.a.c.s.
if and only if $\chi (N)$ is even. This gives one direction. Next
assume $\chi (N)$ is even and let $\xi $ denote a s.a.c.s.. We may
assume $c_1(\xi )=0$ (replace $\xi $ by $\xi + L - \bar
L$, where $L$ is a complex line bundle with $2\cdot c_1(L)=-c_1(\xi
)$). It is well-known that $\xi $ induces an almost complex structure
if and only if $c_4(\xi )$ is equal to the Euler class $e(N)$. We now
use that $N$ admits a topological $Spin(7)$-reduction. In
cohomological terms the existence of such a reduction is equivalent to
$8\cdot e(N)=4\cdot p_2(N)-p_1(N)^2$ (cf. for example \cite{LaMi89},
p. 349). Expressing the Pontrjagin classes of $N$ in terms of the
Chern classes of $\xi $ we get $c_4(\xi )=e(N)$. Thus $\xi $ induces
an almost complex structure on $N$ with vanishing first Chern
class.\proofend

\bigskip
\noindent
{\bf Proof of Proposition \ref{su}:} By Theorem \ref{futaki} the strongly scalar-flat manifold $M$ splits as a Cartesian product of manifolds
with holonomy contained in the special unitary group and manifolds
with $Spin(7)$-holonomy. The former are complex manifolds with
vanishing first Chern class. The latter support almost complex
structures with vanishing first Chern class by Lemma \ref{acs} and the
assumption. For the last statement we recall from \cite{Hi87}, p. 777,
that any $4k$-dimensional almost complex manifold $M$ satisfies the
congruence $sign(M)\equiv (-1)^k\cdot \chi (M)
\bmod 4$.
\proofend

\bigskip
\newpage
\begin{center}4. {\it Relations to the free loop space}\end{center}
\setcounter{section}{4}
\setcounter{theorem}{0}
Some of the previous results indicate that the first Pontrjagin class 
of a manifold $M$ takes a special role in questions concerning 
strongly scalar-flat metrics or Ricci-flat metrics with generic 
holonomy. On the other hand the vanishing of $p_1(M)$ is closely 
related to the existence of a $Spin$-structure on the free loop space 
${\cal L}M$ of $M$.

Motivated by a paper of Stolz (cf. \cite{St95}) we rephrase Theorem 
\ref{firstpont} and Observation \ref{holonomy} in terms of ${\cal L}M$ 
at the end of this section. In \cite{St95} Stolz gives a heuristical 
treatment of a `Weitzenb\"ock formula' for the free loop space of a 
$Spin$-manifold $M$ with $\frac {p_1}2(M)=0$ and conjectures that the 
Witten genus of $M$ vanishes if $M$ admits a metric with positive 
Ricci curvature. The conjecture implies the existence of 
simply-connected Riemannian manifolds with positive scalar curvature 
which do not admit metrics with positive Ricci curvature. 

Let $M$ be a closed oriented connected $n$-dimensional Riemannian 
manifold with free loop space ${\cal L}M=C^\infty (S^1,M)$, where we 
identify $S^1$ with $\R 
/\Z $. Following \cite{St95} we define the scalar curvature of 
${\cal L}M$ to be the map $scal_{{\cal L}M}:{\cal L}M\to \R $, which 
assigns to a loop $\gamma \in {\cal L}M$ the average of the Ricci 
curvature in direction of $\gamma $, i.e. 
$$scal_{{\cal L}M}(\gamma ):=\int_{0}^1Ric (\gamma ^\prime (z))dz.$$
Hence $scal_{{\cal L}M}$ is positive on non-constant loops if $M$ has 
positive Ricci curvature. Also ${\cal L}M$ is scalar-flat, i.e. 
$scal_{{\cal L}M}$ vanishes identically, if $M$ is Ricci-flat. 
Conversely, if $Ric (v)\neq 0$ for some unit tangent vector it is easy 
to find a loop $\gamma $ with $\gamma ^\prime (0)=v$ such that 
$scal_{{\cal L}M}(\gamma)\neq 0$. Hence, if ${\cal L}M$ is scalar-flat 
then $M$ must be Ricci-flat. 

In the following let $M$ be a simply-connected manifold which carries 
a $Spin$-structure given by a $Spin(n)$-principal bundle $P\to M$. The 
free loop space ${\cal L}M$ of $M$ is called $Spin$ if the ${\cal 
L}Spin(n)$-principal bundle ${\cal L}P\to {\cal L}M$ admits a 
reduction to the basic central extension of ${\cal L}Spin(n)$ by $S^1$ 
(for $n\neq 4$ this is the universal extension). As shown in 
\cite{Mc92} the free loop space is $Spin$ if $\frac {p_1}2(M)$ 
vanishes. Here $\frac {p_1}2$ denotes the generator of 
$H^4(BSpin(n);\Z )$ which is half of the universal first Pontrjagin 
class. The converse holds at least if $M$ is two-connected. 

Next assume ${\cal L}M$ is $Spin$. Guided by the finite dimensional 
situation one expects that ${\cal L}M$ carries a `Dirac operator' 
acting on sections of the spinor bundle associated to the 
$Spin$-structure of ${\cal L}M$. In \cite{Wi86} Witten applied 
formally the Lefschetz fixed point formula to this hypothetical Dirac 
operator and derived a bordism invariant of the underlying manifold 
$M$. This invariant, known as the Witten genus, expresses the 
equivariant index of the `Dirac operator' localized at the constant 
loops $M\subset {\cal L}M$. 

In the finite dimensional setting the $\alpha $-invariant, i.e. the 
index of the ordinary Dirac operator, is an obstruction to positive 
scalar curvature metrics on a $Spin$-manifold by the Weitzenb\"ock 
formula (cf. \cite{Hi74}). In analogy one expects that there exists an 
obstruction (namely the Witten genus) to positive scalar curvature on 
the free loop space if ${\cal L}M$ is $Spin$. The first part of the 
following corollary shows that if one replaces positive scalar 
curvature by scalar-flatness the $\alpha $-invariant is an 
obstruction. 

\begin{corollary} Let $M$ be a Riemannian
manifold for which its free loop space ${\cal L}M$ is simply-connected 
and $Spin$.\begin{enumerate} 
\item Assume $M$ has dimension $\geq 5$. If ${\cal L}M$ is scalar-flat then $\alpha (M)=0$, i.e. $M$ admits a metric with positive scalar 
curvature. 
\item Assume $M$ is irreducible. If ${\cal L}M$ is scalar-flat then $M$
is Ricci-flat with generic holonomy group. 
\end{enumerate}
\end{corollary}
\noindent {\bf Proof:} Note that the homotopy groups of ${\cal L}M$ may be computed from the homotopy groups of 
$M$ using the path fibration $\Omega M \to \Gamma M\to M$ and the 
fibration $\Omega M\to {\cal L}M \to M$. In particular, $M$ is 
two-connected since ${\cal L}M$ is simply-connected. As mentioned 
before the existence of a $Spin$-structure for ${\cal L}M$ implies 
that $p_1(M)$ vanishes. Hence the statements follow from Theorem 
\ref{firstpont} and Observation \ref{holonomy}.\proofend 

\begin{finalremark}
\begin{enumerate}
\item Some of the results of this paper have coun\-terparts for non-trivial fundamental groups. 
For example Theorem \ref{firstpont} is still true if $\pi 
_1(M)$ is a finite group for which the Gromov-Lawson-Rosenberg 
conjecture holds (cf. {\rm \cite{De00}}, cf. also the recent paper 
{\rm \cite{BoMc99}} of Botvinnik and McInnes for relations between 
holo\-nomy groups and the Gromov-Lawson-Rosenberg conjecture).
\item The class of simply-connected strongly scalar-flat manifolds of dimension $n\geq 5$
which satisfy an upper diameter bound and a lower curvature bound 
contains only finitely many diffeomorphism classes (for details cf. 
{\rm \cite{De00}}). This follows from the results of 
Cheeger-Fukaya-Gromov on collapsed manifolds in {\rm \cite{ChFuGr92}}. 
\end{enumerate}

\end{finalremark}
\bigskip
\noindent
{\it Acknowledgements:} I thank D. Joyce, J. Lohkamp and U. Semmelmann 
for helpful discussions. Part of this work was completed during a 
visit at the SFB 478 of the University of M\"unster. I am grateful to 
the SFB for its hospitality.

Anand Dessai\\ Department of Mathematics\\University of
Augsburg\\D-86135 Augsburg\\Germany\\ e-mail:
dessai@math.Uni-Augsburg.DE\\ URL:
www.math.Uni-Augsburg.DE/geo/dessai/homepage.html

\end{document}